\documentclass[12pt,a4paper]{article}
\usepackage{graphicx}
\begin{document}
 
 \thispagestyle{empty}
 
 \title{Platonic polyhedra tune the 3-sphere:\\ Harmonic analysis on simplices.}
 \author{Peter Kramer, Institut fuer Theoretische Physik,\\ U Tuebingen, Germany.}
 \maketitle

\section*{Abstract.}
A spherical topological manifold of dimension $n-1$ forms a prototile  on its  cover,
the (n-1)-sphere. The tiling is generated by the fixpoint-free 
action of the group of deck transformations. By 
a general theorem, this group is isomorphic to the first homotopy group.  
A basis for the harmonic analysis on the (n-1)-sphere  is given by the spherical harmonics
which transform according to irreducible representations of the orthogonal group.
Multiplicity and selection rules appear in the form of
reduction of group representations.
The deck transformations form a  subgroup, and so the representations 
of the orthogonal group  can be reduced to those of this subgroup. Upon reducing to the identity representation  of the subgroup, the reduced subset of spherical harmonics becomes periodic on the tiling and 
tunes the harmonic analysis on the (n-1)-sphere  to the manifold. A particular class of spherical 3-manifolds arises
from the Platonic polyhedra. The harmonic analysis on the Poincare dodecahedral 3-manifold 
was analyzed along these lines. For comparison we construct here the harmonic analysis on simplicial spherical manifolds of dimension $n=1,2,3$. 
Harmonic analysis can be applied to the 
cosmic microwave background observed in astrophysics.
Selection rules found in this analysis can detect   
the multiple connectivity  of spherical 3-manifolds on the space part of cosmic space-time.

\vspace{0.2cm}

\section{Introduction.}

Viewed on its universal cover $S^{(n-1]}$, a spherical topological manifold ${\cal M}$ is the prototile of
a tiling generated by the corresponding deck transformations. By a general theorem in topology given 
by Seifert and Threlfall \cite{SE34} pp. 195-8, the group $deck({\cal M})$ of deck transformations is
isomorphic to the first homotopy group $\pi_1({\cal M})$ and hence is a topological invariant.

A particular class of spherical  3-manifolds arises  from the five Platonic polyhedra. 
Everitt \cite{EV04} discusses their topology  and gives a graphical algorithm for their homotopy groups.
In Table 1.1 we list the five polyhedra and from \cite{EV04} give the known order of their homotopy group.

{\bf Table 1.1}. The five 3-manifolds arrising from the Platonic polyhedra, the orders 
$|\pi_1({\cal M}|$ of their first homotopy groups, and the volume fraction 
$frac({\cal M})= |\pi_1({\cal M})|^{-1}$ of the prototile 
w.r.t the volume of the 3-sphere.

\begin{equation}
 \begin{array}{llllll}
{\cal M}:        &{\rm tetrahedron}&{\rm cube}&{\rm octahedron}&{\rm icosahedron}&{\rm dodecahedron}\\
|\pi_1({\cal M})|:&5          &8   &8         &&120\\ 
frac({\cal M}):&0.2      &0.125&0.125    &&0.0083..
 \end{array}
\end{equation}

\vspace{0.3cm}

Since the groups $deck({\cal M})$ and $\pi_1({\cal M})$ are isomorphic and $deck({\cal M})$
acts fixpoint-free on $S^3$, the volume fraction of ${\cal M}$ taken as prototile 
on $S^3$ is equal to $|\pi_1({\cal M})|^{-1}$. We see that the tetrahedron and the dodecahedron 
display extremal values of this fraction.

For general notions of topology we refer to \cite{SE34}, \cite{TH97}.
The harmonic analysis on these manifolds can be started from $S^{(n-1)}$.
There its basis is the complete, orthonormal set $\langle Y^{\lambda}\rangle$ of spherical harmonics, the square integrable 
eigenmodes of $S^{(n-1)}$. To pass to a 3-manifold ${\cal M}$ universally covered by $S^{(n-1)}$,
consider the maximal subset $\langle Y^{\lambda 0}\rangle$ of this basis periodic with respect to deck transformations. Due to the periodicity, it can be restricted 
to the prototile ${\cal M}$ and forms its eigenmodes. These periodic  eigenmodes tune the sphere $S^{(n-1]}$ to the topology of
${\cal M}$.

To analyze in detail the periodic eigenmodes we turn to groups and their
representations. Under the group $O(n,R)$ of isometries of $S^{(n-1)}$,
the  spherical harmonics $ Y^{\lambda} $ transform according to a set $D^{\lambda}$ of
irreducible representations . The periodic subset $Y^{\lambda 0}$ transforms 
according to the identity representation $D^0$ of the group $deck({\cal M})$ 
of deck transformations. In terms of the group/subgroup pair $O(n,R)>deck({\cal M})$,
we require the reduction of the irreducible representations $D^{\lambda}>D^0$.  Clearly not all representations $D^{\lambda}$ will reduce to the 
representation $D^0$. The non-occurrence provides  a selection rule 
on eigenmodes of $S^{(n-1)}$ which is another mark for the topology of  ${\cal M}$. 

Motivated by physics, harmonic analysis on topological 3-manifolds  has  been invoked in cosmological models
of the space part of space-time \cite{LA95}, \cite{LU03}. A direct experimental access 
to the topology from the autocorrelation of the cosmic matter distribution is difficult.
As an alternative, the rich data from fluctuations of the 
Cosmic Microwave Background (CMB) radiation were examined by harmonic analysis. It is hoped to find in this way the characteristic selection rules and tuning for a specific non-trivial topology, distinct from the standard simply-connected one. Of the Platonic 3-manifolds,
the Poincare dodecahedral manifold of minimal volume fraction and its eigenmodes have  found particular attention. Representation theory was applied to the harmonic analysis on Poincare's dodecahedral 3-manifold in \cite{KR05}. A comparative study of the harmonic analysis, tuned to different  topological 3-manifolds, can provide clues for future applications.

To initiate such a comparative study,  we turn here to simplicial manifolds on $S^{(n-1)},\: n-1=1,2,3$.
For illustration of the group and representation theory we start in sections 2 and 3 with the cases  $n-1=1,2$.  Section 4 deals with the Platonic tetrahedral 3-manifold.

A regular simplex with $(n+1)$ vertices can be centrally projected to $S^{(n-1)}$ to yield a tiling 
into $(n+1)$ spherical simplices. The simplicial tiling can be  generated
by the fixpoint-free action of the cyclic group $C_{n+1}$ acting as 
${\rm deck}({\cal S}_0(n-1))$ on a prototile. In Fig. 1 we illustrate symbolically the tilings and simplicial manifolds
for $n-1=1,2,3$. 

We require
the identity representation $D^0$ of $C_{n+1}$ in the reduction of representations 
of the groups
\begin{equation}
\label{t1}
 O(n,R) > C_{n+1}.
\end{equation}
Working with this pair of groups will display the difference in the topology 
between $S^{(n-1)}$ and ${\cal S}_0(n-1)$ as part of the harmonic analysis.

\begin{center}
\includegraphics{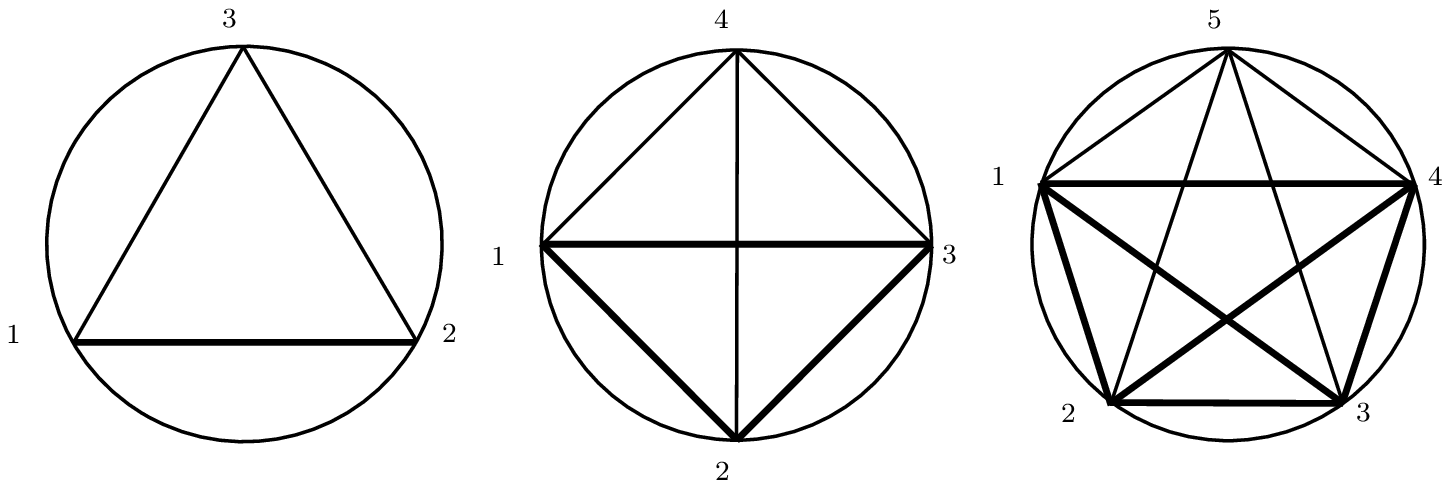} 
\end{center}

Fig.1 Simplicial topological manifolds on spheres. Regular  simplices with $(n+1)$ vertices are inscribed into
spheres $S^{(n-1)}, n=2,3,4$, here  represented symbolically  by a circle. Central  projection yields a simplicial tiling of $S^{(n-1)}$.
Projection of the subsimplex ${\cal M}={\cal S}_0(n-1)$, heavy 
lines, obtained by dropping the vertex labelled $(n+1)$,  represents a simplicial topological $(n-1)$-manifold as a prototile on its universal cover $S^{(n-1)}$.
The group $deck({\cal S}_0(n-1))$ is the cyclic group $C_{n+1}$. It acts fixpoint-free 
on $S^{(n-1)}$ and produces the simplicial tiling. By a general theorem it is isomorphic to the first homotopy group
$\pi_1({\cal S}_0(n-1))$ which is a topological invariant. Restricting the basis of the harmonic analysis from $S^{(n-1)}$ to its subset periodic under $C_{n+1}$ we construct the harmonic analysis on ${\cal S}_0(n-1)$. 
\vspace{0.3cm}

Now we note that the original $n$-simplex, of which ${\cal S}_0(n-1)$ forms a spherical face,
has the symmetry group $S(n+1)$, the symmetric group on $(n+1)$ variables. The representation theory of $S(n+1)$ is well known \cite{RO61} pp. 38-9, \cite{HA62} pp. 214-31, and we shall
make full use of it.  

$S(n+1)$ is 
a Coxeter group generated by reflections, see \cite{HU90}.
It has the Coxeter diagram with $n$ nodes 
\begin{equation}
\label{t3}
S(n+1) = \circ-\circ-\ldots -\circ 
\end{equation}
Each node of this diagram, \cite{HU90} pp. 31-3, describes the generator of a Weyl reflection $W_a$ determined 
by a Weyl or reflection vector $a$ in the Euclidean space $E^n$ embedding $S^{(n-1)}$, with $W_a$ the involutive map
\begin{equation}
\label{t3a}
x \in E^n: (W_a \times x)  \rightarrow  
x-2\frac{\langle x,a\rangle}{\langle a,a\rangle}\:  a. 
\end{equation}
Moreover the diagram eq. \ref{t3} implies definite relations and scalar products between 
the various Weyl reflection vectors. In section 4.3 we shall take
$S(5)$ as a Coxeter group and describe its Weyl reflection vectors.

. 

Since $C_{n+1}$
is generated by the cyclic permutation $(1,2,\ldots, n+1) \in  S(n+1)$, we can refine the group/subgroup 
pair eq. \ref{t1} as 
\begin{equation}
\label{t2}
O(n,R) > S(n+1) > C_{n+1}. 
\end{equation}
The group $S(n+1)$ here is   taken as a subgroup of $O(n,R)$ via the orthogonal irreducible 
representation $D^{[n1]}$ of partition $f=[n1]$ since
(i) the Weyl reflections of the defining representation  for the Coxeter group eq. \ref{t3} are orthogonal transformations, (ii) the irreducible orthogonal representations of $S(n+1)$ are characterized by Young diagrams $f$. The orthogonal irreducible representation $D^{f}, f=[n1]$ has dimension $n$ and is 
equivalent to the defining representation of the Coxeter group.

The cyclic subgroup $C_{n+1} < S(n+1)$ in eq. \ref{t2} is generated by the product of
all the Weyl reflection generators of the Coxeter group. This product is termed the Coxeter element, see \cite{HU90} pp. 74, 174.

An advantage in using the scheme eq. \ref{t2} is the fact that the subgroup 
$S(n)<S(n+1)$ permuting the
vertices $1,2,\ldots, n$ naturally appears as a symmetry subgroup acting on the spherical simplex ${\cal S}_0(n)$. Borrowing  the terminology from space groups in $E^n$,
we call $S(n)$ the point group of ${\cal S}_0(n)$. 

The reduction of irreducible representations for $O(n,R)>S(n+1)$ was studied in \cite{KM66} in full detail for $n=2,3$. This allows to work out the harmonic  analysis for these cases. 
We apply the representation theory of groups, following \cite{KM66}, in the following sections  to the harmonic analysis on the  simplicial manifolds for $n-1=1,2,3$ .
By a first  step in the reduction of representations according to  eq. \ref{t2} we find multiplicities and selection rules 
for the representations of $O(n,R)$ and $S(n+1)$ which subduce to the identity representation 
of $C_{n+1}$. This reduction we handle by character technique. 

In the second step we  characterize and construct an explicit set of 
orthonormal eigenmodes which span a unitary linear space of functions ${\cal L}^2$  for the harmonic analysis on the simplicial topological manifolds. The measure $d\mu$ on ${\cal L}^2$
in our analysis is taken over from the universal cover $S^{(n-1)}$. There are two alternative 
but equivalent approaches differing in the choice of the domain:

(I:) A given function  ${\cal F}$ with domain ${\cal S}_0(n-1)$ can be extended 
to a $C_{n+1}$-periodic function on the universal cover $S^{(n-1)}$ and then analyzed exclusively in terms of a $C_{n+1}$-periodic 
basis. 

(II:) By $C_{n+1}$-periodicity, this  analysis can be  collapsed  entirely to the domain ${\cal S}_0(n-1)$, which forms a spherical tile of $S^{(n-1)}$ of volume ${\rm vol}(S^{(n-1)})/(n+1)$,
with the range of integration 
restricted to the topological manifold as domain .  
The expansion coefficients of a given function ${\cal F}$ with domain ${\cal S}_0(n)$ 
are found from  the scalar products on ${\cal L}^2$ of ${\cal F}$ with the basis functions. 

We focus on the choice (I) since the modification to (II) is straightforward.
The main results of the following three sections are given in Tables after each section.

\section{The 1-simplex ${\cal S}_0(1)$ on the circle $S^1$.}
As the simplest paradigm we treat here the spherical simplex ${\cal S}_0(1)$,  with universal
cover the unit circle $S^1$ according to Fig. 1, left.

\subsection{The symmetric group $S(3)$.}
Using orthogonal coordinates in the Euclidean plane $E^2$, the symmetry group $S(3)$ of the 
regular triangle 
is generated by two Weyl reflections with matrices
\begin{equation}
\label{t4}
(1,2):
\left[
\begin{array}{ll}
-1&0\\
0&1\\ 
\end{array}
\right],\;
(2,3):
\left[
\begin{array}{ll}
-\frac{1}{2}&\frac{1}{2}\sqrt{3}\\ 
\frac{1}{2}\sqrt{3}&\frac{1}{2}
\end{array}
\right]
\end{equation}
These matrices generate  Young's orthogonal representation $D^{[21]}$ with basis vectors corresponding  to the two Young tableaus
\begin{equation}
\label{t4a} 
\left[
 \begin{array}{ll}
1&3\\
2&\\
\end{array}  
\right],\:
\left[
 \begin{array}{ll}
1&2\\
3&\\
\end{array}  
\right].
\end{equation}

The two other irreducible representations of $S(3)$ are
\begin{eqnarray}
\label{t5}
&&D^{[3]}:  D^{[3]}(1,2)=D^{[3]}(2,3)=1,
\\ \nonumber 
&&D^{[111]}:  D^{[111]}(1,2)=D^{[111]}(2,3)=-1,
\end{eqnarray}
The character table of $S(3)$ in terms of classes, denoted by the cycle structures
$k= (1)^3, (2)(1), (3)$ is given as {\bf Table 2.1}. 

\subsection{The cyclic group $C_3$.}
The cyclic group as a subgroup $C_3<S(3)$ is given  in cycle notation by 
\begin{equation}
\label{t7}
 C_3= \{(1,2,3),(1,2,3)^2=(3,2,1),(1,2,3)^3=e\}.
\end{equation}
In terms of the complex number $\lambda_3=\exp(\frac{2\pi i}{3})$, the 
three irreducible representations $D^{\alpha}, \alpha=0,1,2$ of $C_3$ coincide with their  characters and are 
given for the group elements by {\bf Table 2.2}.

\subsection{The reduction $S(3)>C_3$.}
The multiplicity $m(f,\alpha)$ of the subgroup irreducible representation $D^{\alpha}$ in the
group irreducible representation $D^f$ for $H<G$ is given from character technique by
\begin{equation}
\label{t9}
m(f,\alpha)= \frac{1}{|H|}\sum_{k\in H} n(k)\: \chi^f(k)\overline{\chi}^{\alpha}(k).
\end{equation}
Here $n(k)$ is the number of elements in class $k$ of the subgroup $H$.
We are interested only in the identity representation $D^0$ of $H=C_3$,  
collect the characters $\chi^f(k)$ of the group elements eq. \ref{t7} in the irreducible representations 
of $S(3)$ and  $\chi^0(k)=D^0(k)=1$, and compute the multiplicities $m(f,0)$ from eq. \ref{t9} with $|C_3|=3,\: n(k)=1$ in  {\bf Table 2.3}.

As a result, the identity representation $D^0$ of $C_3$ is contained once in  the irreducible representations $D^{[3]}, D^{[111]}$ of $S(3)$ but not in $D^{[21]}$.

\subsection{The reduction $O(2,R)>S(3)$.}
Now we turn to the group/subgroup pair $O(2)>S(3)$. The harmonic analysis on the 
circle $S^1$ is given by the complex Fourier series, spanned by the complete orthonormal 
set  of functions of $z=\exp(i \phi),\: 0\leq \phi <2\pi$,
\begin{eqnarray}
 \label{t11}
&&Y_m(z)=\frac{1}{\sqrt{2\pi}}z^m, m=0, \pm 1,\; \ldots,
\\ \nonumber
&&\int_0^{2\pi} d\phi \overline{Y_p}(z) Y_q(z)= \delta_{pq},
\\ \nonumber
&& T_{(1,2)}Y_m(z)= (-1)^m Y_{-m}.
\end{eqnarray}
We include the reflection $(1,2): \phi \rightarrow \pi-\phi$ by introducing the new 
basis 
\begin{equation}
\label{t12}
Y_0,\: Y_{m,\epsilon}= \sqrt{\frac{1}{2}}\left[Y_m+\epsilon (-1)^m Y_{-m}\right],\;
\epsilon= \pm 1,\; m=1,2,\ldots
\end{equation}
We are now ready to give the reduction of representations in $O(2)>S(3)$ from \cite{KM66}.
We introduce $\nu$ by
\begin{equation}
\label{t13}
\nu:\; m  \equiv \nu\; mod\; 3,\; \nu=0,1,2.
\end{equation}
Then the results from \cite{KM66} pp. 255-7 can be combined with those from 
{\bf Table 2.3} to yield the full reduction of representations for the groups 
eq. \ref{t2} in {\bf Table 2.4}.

{\bf Table 2.4} determines the harmonic analysis on ${\cal S}_0(1)$ carried out on $S^1$:
Any function ${\cal F}$ with domain ${\cal S}_0(1)$ has a unique 
$C_3$-periodic extension to $S^1$ and can be expanded into $C_3$-periodic
basis functions. The  alternative would be  to  restrict
the  integration in eq. \ref{t11} to the sector ${\cal S}_0$ of angular range $2\pi/3$ and 
to reduce the volume from $2\pi$ to $2\pi/3$. 

The $C_3$-periodic basis is obtained from  the  irreducible representations of $O(2,R)$ by restriction to the values 
$\nu=0,\; m \equiv 0\; mod\; 3, \epsilon=\pm 1$. This restricts the irreducible representations of $S(3)$ 
to $f=\left[3\right], \left[111\right]$. Finally the representation of the point symmetry group $S(2)$ of 
the spherical simplex ${\cal S}_0(1)$ is given in terms of partitions $f$ by 
\begin{equation}
\label{t15}
m=0: \left[2\right],\; m>0, \epsilon=1: f=\left[2\right],\; \epsilon=-1: f=\left[11\right]. 
\end{equation}

\subsection{Tables for $O(2,R)>S(3)>C_3$.}

{\bf Table 2.1}. Characters for the $3$ partitions $f$ and classes $k$ in cycle notation 
of the symmetric group $S(3)$.

 \begin{equation}
\label{t6}
\begin{array}{r|rrr}
\chi^f(k) & (1)^3&(2)(1)&(3)\\
\cline{1-4}
\left[3\right]       &1&1&1\\
\left[21\right]      &2&0&-1\\
\left[111\right]     &1&-1&1
\end{array}.
 \end{equation}
\vspace{0.3cm}

{\bf Table 2.2}. Characters for the $3$ irreducible representations of $C_3$ with $\lambda_3=\exp(\frac{2\pi i}{3})$.
\begin{equation}
\label{t8}
\begin{array}{l|lll}
\chi^{\alpha} &e&(1,2,3)&(3,2,1) \\
\cline{1-4}
D^0  & 1&1&1\\
D^1  & 1&\lambda_3&\lambda^2_3\\
D^2  & 1& \lambda^2_3&\lambda_3\\
\end{array}
\end{equation} 
\vspace{0.3cm}

{\bf Table 2.3}. Characters and multiplicity $m(f,0)$  in the reduction of 
representations $S(3) >C_3$.

\begin{equation}
\begin{array}{l|lll|l}
k & e& (1,2,3) & (3,2,1)& m(f,0)\\
\hline
\chi^{\left[3\right]}(k)& 1 & 1& 1& 1\\
\chi^{\left[21\right]}(k)& 2 & -1& -1& 0\\
\chi^{\left[111\right]}(k)& 1 & 1& 1& 1\\
\chi^0(k)&1&1&1&
\end{array}
\end{equation}
\vspace{0.3cm}

{\bf Table 2.4}. Reduction of irreducible representations for the groups 
$O(2,R)>S(3)>C_3$.

\begin{equation}
 \label{t14}
\begin{array}{ll|l}
O(2)& S(3) & C_3\\
(m,\epsilon)& f & m(f,0)\\
\cline{1-3}
m=0 & \left[3\right] & 1\\
m>0, \nu=0, \epsilon=1& \left[3\right]&1\\
m>0, \nu=0, \epsilon=-1& \left[111\right]&1\\
m>0, \nu=1, \epsilon=\pm 1& \left[21\right]&0 \\
m>0, \nu=2, \epsilon=\pm 1& \left[21\right]&0 \\
\end{array}
\end{equation}
\vspace{0.3cm}

\section{The 2-simplex ${\cal S}_0(2)$ on the sphere $S^2$.}

The 2-faces of the regular 3-simplex centrally projected to the sphere $S^2$, Fig. 2, tiles it 
into 4 spherical triangles. The tiling has the symmetry $S(4)$ with the Coxeter diagram
\begin{equation}
 \label{t16}
S(4): \circ-\circ-\circ.
\end{equation}

\begin{center}
\includegraphics{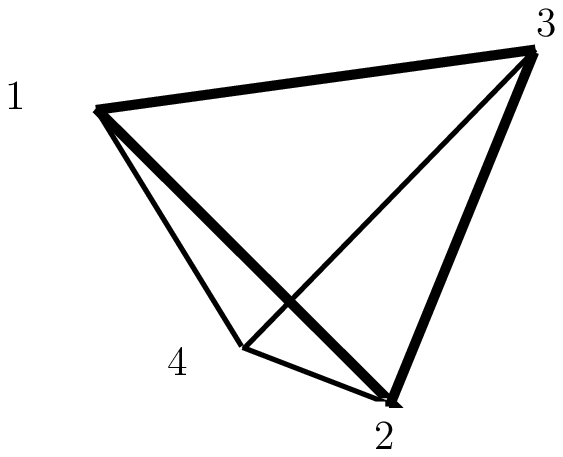} 
\end{center}

Fig.2 The 3-simplex with  vertices $1,2,3,4$ can be centrally projected to the sphere $S^2$. The spherical 
triangle  spanned by the vertices $1,2,3$ on $S^2$ is taken as the simplicial topological manifold 
${\cal S}_0(2)$.
\vspace{0.4cm}

The triangle obtained by dropping the vertex labelled 4
we take as the simplicial topological manifold ${\cal S}_0(2)$ seen on its universal cover 
$S^2$. The generators $(1,2),(2,3),(3,4)$ of $S(4)$ act as three reflection planes of the 
tetrahedron. The first two of them map ${\cal S}_0(2)$ into itself and generate the 
point group $S(3)$ of the simplicial manifold. 

\subsection{The group $S(4)$.}

The characters of the irreducible representations 
of $S(4)$ are given from \cite{HA62} p. 187 in {\bf Table 3.1}.

\subsection{The cyclic group $C_4$.}
The cyclic group $C_4$ is generated as 
\begin{equation}
\label{t18}
C_4: \left\{ (1,2,3,4), (1,2,3,4)^2=(1,3)(2,4),(1,2,3,4)^3=(4,3,2,1), (1,2,3,4)^4=e\right\} 
\end{equation}
From these expressions it can be seen that the elements of $C_4$ generate all the cosets 
of the subgroup $ S(3)<S(4)$. Hence their actions produce the full tiling of $S^2$ 
by triangles, are fixpoint-free and so fulfill all the properties required for
$deck(\tilde{{\cal S}}_0(2))$.

With $\lambda_4=\exp(2\pi i/4)=i$, the 4 $1$-dimensional representations of $C_4$ are generated by
\begin{equation}
\label{t19}
D^0(1,2,3,4)=1,  D^1(1,2,3,4)=i,D^2(1,2,3,4)=i^2=-1,D^3(1,2,3,4)=i^3=-i,
\end{equation}
The characters of $C_4$ are given in {\bf Table 3.2}.

\subsection{The reduction $S(4)>C_4$.}
We require the reduction of the representations $f$ of $S(4)$ to the identity 
representation $D^0$ of $C_4$. For this purpose we write the characters $\chi^f$ of these 
representations and $\chi^0=D^0$ for the classes $k$ of $C_4$ in the form of {\bf Table 3.3} and
compute  the multiplicity $m(f,0)$.

As a result, the identity representation $D^0$ of $C_4$ is contained once in the representations 
$f=\left[4\right],\left[22\right], \left[211\right]$ but not 
in $f=\left[31\right],\left[1111\right]$. Since the representations $f=\left[22\right], \left[211\right]$ are $2$- and $3$-dimensional, we must still find the 1-dimensional subspaces 
for the representation $D^0(C_4)$. This will be done by projection in subsection 3.5.1, 3.5.2.

\subsection{The reduction $O(3,R)>S(4)$.}
Next we turn to the harmonic analysis of $O(3,R)$ acting on the sphere $S^2$. 
$O(3)$ is the direct product $SO(3,R) \times {\cal J}$, where ${\cal J}$ is 
the group generated by parity. Following \cite{KM66}, we characterize 
the irreducible representations by $(l, \kappa)$ where $l=0,1,\ldots $
is the integer angular momentum  and $\kappa=\pm, 1$ 
a parity label. The basis of the irreducible representation $D^l(SO(3,R))$
is spanned in spherical polar coordinates $\theta, \phi$ by the spherical harmonics, 
\cite{ED57} pp.19-25,
\begin{eqnarray}
 \label{t20a}
&&Y_{lm}(\theta,\phi),\: -l \leq m \leq l,
\\ \nonumber
&&\int \overline{Y}_{l'm'} Y_{lm} \sin(\theta) d\theta d\phi= \delta_{l'l} \delta_{m'm}.
\end{eqnarray}

which under 
the parity operation $P$ transform as 
\begin{equation}
\label{t21}
P Y_{lm}= (-1)^l Y_{lm},\:   -l\leq m \leq l.
\end{equation}
The parity therefore is
\begin{equation}
 \label{t22}
\kappa= (-1)^l. 
\end{equation}

In {\bf Table 3.4}  we reproduce from \cite{KM66} pp. 259-60
the decomposition of irreducible representations for $(O(3,R)\times {\cal J})>S(4)$ and 
combine it with the results from {\bf Table 3.3}. Because of eq. \ref{t22},
parity $\kappa$ and $l$ are correlated. 

\vspace{0.4cm}
This Table can be extended by use of character technique given in  \cite{KM66}.

Counting the total number of states up to $l=4$ from {\bf Table 3.4} one finds $25$ basis states 
for the group $O(3,R)$. Only $7$ of them contribute to the $C_4$-periodic subset. 

\subsection{The explicit  reduction $O(3,R)>S(4)>C_4$.}
For higher representations $(l,\kappa)$ the multiplicity $m((l, \kappa),f)$  of the partition
$f$ will take values larger than $1$. In \cite{KM66} we devised  a way which allows
for a complete orthogonal basis in this case. The idea is to introduce an additional
hermitian operator, termed a generalized Casimir operator, as a polynomial in the components of the angular momentum operator,
which for fixed $(l,\kappa)$ by distinct eigenvalues separates repeated partitions $f$.
We refer to \cite{KM66} pp. 263-66 for the details.

We now  determine the $1$-dimensional subspaces for the partitions 
$f=\left[22\right], \left[211\right]$ belonging to the representation $D^0$
of $C_4$.

\subsubsection{The partition $f=\left[211\right]$.}

Following \cite{KM66}, eq. (6.9), we choose particular coordinates with respect to the 4-simplex Fig. 2,
\begin{eqnarray}
 \label{t24}
&&y_1=\frac{1}{2}\left[x_1+x_2-x_3-x_4\right]
\\ \nonumber
&&y_2=\frac{1}{2}\left[x_1+x_3-x_2-x_4\right]
\\ \nonumber
&&y_3=\frac{1}{2}\left[x_1+x_4-x_2-x_3\right]
\end{eqnarray}
These coordinates are convenient and stem from the analysis of the
tetrahedral group $T_d$ isomorphic to $S(4)$ \cite{KM66} p. 258-68. 

Expressed in these coordinates, the generators of $S(4)$ have the $3 \times 3$ matrix representation
\begin{equation}
 \label{t25}
D^{\left[31\right]'}(1,2)= 
\left[
\begin{array}{lll}
1&0&0\\
0&0&-1\\
0&-1&0
\end{array}
\right]\:
D^{\left[31\right]'}(2,3)= 
\left[
\begin{array}{lll}
0&1&0\\
1&0&0\\
0&0&1
\end{array}
\right]\:
 D^{\left[31\right]'}(3,4)= 
\left[
\begin{array}{lll}
1&0&0\\
0&0&1\\
0&1&0
\end{array}
\right].
\end{equation}
which generate the representation $D^{\left[31\right]'}$ used in relation with 
the isomorphic tetrahedral group \cite{KM66} p. 261. We use primes to distinguish 
this representation from the equivalent Young representation $D^{\left[31\right]}$ \cite{HA62} pp. 225-6. 

The representation 
with the associate partition
$D^{\left[211\right]'}$ 
is obtained by multiplying each matrix for a generator in eq. \ref{t25} by $(-1)$, corresponding to the 
representation $D^{\left[1111\right]}$.
Within this representation we construct the matrix of the Coxeter element
\begin{equation}
 \label{t26}
D^{\left[211\right]'}(1,2,3,4)= 
 \left[
\begin{array}{lll}
0&0&-1\\
0&1&0\\
1&0&0
\end{array}
\right].
\end{equation}
and by taking powers the representation for all $4$ elements 
$h\in C_4$. The  projection operator to the representation $D^0$ of $C_4$ becomes
\begin{equation}
\label{t27} 
P^{\left[211\right]',0}=\frac{1}{4}\sum_{h\in C_4}D^{\left[211\right]'}(h)=
 \left[
\begin{array}{lll}
0&0&0\\
0&1&0\\
0&0&0
\end{array}
\right].
\end{equation}
It picks the second basis function of $D^{\left[ 211\right]'}$ as the one which transforms
according to $D^0$. 

Alternatively we describe the same representation $f=\left[ 211\right]$ in terms of the standard Young-Yamanouchi  
basis as given by Hamermesh \cite{HA62} pp. 224-226. In contrast to the 
primed representation eq. \ref{t25} 
we label the basis vectors of this unprimed representation by three Yamanouchi 
symbols  corresponding to three Young tableaus \cite{HA62} pp 225-6: 
\begin{equation}
\label{t27a}
1: 3211, \left[
 \begin{array}{ll}
1&2\\
3&\\
4&\\
\end{array}  
\right],
2:3121, \left[
 \begin{array}{ll}
1&3\\
2&\\
4&\\
\end{array}  
\right],
3:1321, \left[
 \begin{array}{ll}
1&4\\
2&\\
3&\\
\end{array}  
\right].
\end{equation}

In terms of the canonical group/subgroup 
sequence  $S(4)>S(3)>S(2)$ underlying the Young representation,
the three Young basis vectors in eq. \ref{t27a} correspond to the partitions 
$\left[21\right],\left[21\right],\left[111\right]$ of $S(3)$ in that order.
We find the matrix representations 
of the reflection generators and of the generator of $C_4$ as
\begin{eqnarray}
 \label{t25a}
&&D^{\left[211\right]}(1,2): 
\left[
\begin{array}{lll}
1&0&0\\
0&-1&0\\
0&0&-1
\end{array}
\right]\:
D^{\left[211\right]}(2,3): 
\left[
\begin{array}{lll}
-\frac{1}{2}&\frac{1}{2}\sqrt{3}&0\\
\frac{1}{2}\sqrt{3}&\frac{1}{2}&0\\
0&0&-1
\end{array}
\right]\:
\\ \nonumber 
&& 
D^{\left[211\right]}(3,4): 
\left[
\begin{array}{lll}
-1&0&0\\
0&-\frac{1}{3}&\sqrt{\frac{8}{9}}\\
0&\sqrt{\frac{8}{9}}&\frac{1}{3} 
\end{array}
\right],
D^{\left[211\right]}(1,2,3,4): 
\left[
\begin{array}{lll}
\frac{1}{2}&-\frac{1}{6}\sqrt{3}&\sqrt{\frac{2}{3}}\\
\frac{1}{2}\sqrt{3}&\frac{1}{6}&-\frac{1}{3}\sqrt{2}\\
0&\frac{2}{3}\sqrt{2}&\frac{1}{3}
\end{array}
\right].
\end{eqnarray}
We can find the linear combination of the basis states belonging to $D^0(C_4)$ 
by constructing the eigenvector $\psi^{\left[211\right] 0}$ of the representation matrix $D^{\left[211\right]}(1,2,3,4)$ with eigenvalue $1$. This eigenvector is found 
in terms of the Young tableau labels eq. \ref{t27a} as
\begin{equation}
 \label{t25b}
\psi^{\left[211\right],0}= 
(\sqrt{\frac{1}{2}}, \sqrt{\frac{1}{6}}, \sqrt{\frac{1}{3}}).
\end{equation}
From the  linear combination of the basis states with the coefficients eq. \ref{t25b} we can draw a surprising   conclusion. 
The Young representation basis eq. \ref{t27a}  displays explicitly the sub-partition of
$S(3)$ as the part of the Young diagram occupied by the numbers $1,2,3$.
The basis vector eq. \ref{t25b} of $D^0(C_4)$ in the Young representation
$\left[ 211\right]$ of $S(4)$ is a generic mixture of the two inequivalent  representations 
$\left[ 21\right], \left[ 111\right]$
of the group $S(3)$. Although   we interpreted $S(3)$ as the point group of the simplicial manifold.
the reduction to $D^0$ does not  preserve a single  representation of this subgroup. 
The cyclic group 
$C_4$ still provides four coset generators for the subgroup $S(3)<S(4)$, and so  
$S(4)$ is generated from all the products of  elements from its subgroups
$S(3), C_4$. Moreover our result 
eq. \ref{t25b} shows that the use of the Young representation may not
optimal for the construction.
Comparison of eq. \ref{t25b} with the simpler result eq. \ref{t27}  favours the 
primed tetrahedral form eq. \ref{t25} of the representation.

\subsubsection{The partition $f=\left[22\right]$.}
The second representation $D^{\left[22\right]}$ we treat in the Young orthogonal 
representation \cite{HA62} p. 226 with the Young tableaus
\begin{equation}
\label{t25c}
1:\:
\left[ 
\begin{array}{ll}
1& 2\\
3& 4
\end{array}
\right],\; 
2:\: 
\left[\begin{array}{ll}
1& 3\\
2& 4
\end{array}
\right], 
\end{equation}

and with generators 
\begin{eqnarray}
 \label{t28}
&&D^{\left[ 22\right]}(1,2): 
\left[ 
\begin{array}{ll}
 1 &0 \\
 0 &-1  \\
\end{array}
\right],\:
D^{\left[ 22\right]}(2,3): 
\left[ 
\begin{array}{ll}
 -\frac{1}{2} &-\frac{1}{2}\sqrt{3} \\
-\frac{1}{2}\sqrt{3}  &\frac{1}{2}  \\
\end{array}
\right],\:
\\ \nonumber
&&D^{\left[ 22\right]}(3,4):
\left[ 
\begin{array}{ll}
 1 &0 \\
 0 &-1  \\
\end{array}
\right],\:
D^{\left[ 22\right]}(1,2,3,4):
\left[ 
\begin{array}{ll}
 -\frac{1}{2} &\frac{1}{2}\sqrt{3} \\
\frac{1}{2}\sqrt{3}  &\frac{1}{2}  \\
\end{array}
\right].\:
\end{eqnarray}
By use of eq. \ref{t18} we find from the second line of eq. \ref{t28} the  representation of the Coxeter element. 
The projection operator for the  representation $D^0(C_4)$ becomes the matrix
\begin{equation}
 \label{t29}
P^{\left[22\right],0}=\frac{1}{4}\sum_{h\in C_4} D^{\left[ 22\right]}(h)=
\frac{1}{4}\left[
\begin{array}{ll}
1&\sqrt{3}\\
\sqrt{3}&3\\
\end{array}
\right].
\end{equation}
and determines a unique linear combination in the Young tableaus,
\begin{equation}
\label{t30}
\psi^{\left[22\right],0}= (\frac{1}{2},\sqrt{\frac{3}{4}})
\end{equation}
as the basis for $D^0(C_4)$. The representation of the 
point group $S(3)$ of the simplicial manifold is fixed  by the partition $f=\left[21\right]$.

This completes the group/subgroup basis construction of the harmonic analysis for 
a function ${\cal F}$ on the simplicial
topological manifold ${\cal S}_0(2)$ when  $C_4$-periodically  extended to $S^2$. The new feature here is the 
appearance  of multiplicity problems which require projections and
the construction of appropriate operators. The $C_4$-periodic basis for the harmonic analysis
of ${\cal S}_0(2)$ on the sphere $S^2$ is characterized by the labels $((l,\kappa),f)$ of the irreducible
representations of $O(3,R),S(4)$, with $f$  restricted to allowed values.

\vspace{0.3cm}

\subsection{Tables for $O(3,R)>S(4)>C_4$.}

{\bf Table 3.1}. Characters for the group $S(4)$.

\begin{equation}
\label{t17}
\begin{array}{r|rrrrr}
\chi^f(k) & (1)^4&(4)&(3)(1)&(2)^2&(2)(1)^2\\
\cline{1-6}
\left[4\right]       &1&1&1&1&1\\
\left[31\right]      &3&-1&0&-1&1\\
\left[22\right]      &2&0&-1&2&0\\
\left[211\right]     &3&1&0&-1&1\\
\left[1111\right]    &1&-1&1&1&-1\\
\end{array}.
 \end{equation}

\vspace{0.4cm}

{\bf Table 3.2}. Characters for the cyclic group $C_4$.

\begin{equation}
\label{t19a}
\begin{array}{r|rrrr}
\chi^{\alpha} &e&(1,2,3,4)&(1,3)(2,4)& (4,3,2,1) \\
\cline{1-5}
D^0  & 1&1&1&1\\
D^1  &1 &i &-1&-i\\
D^2  & 1&-1&1&-1\\
D^3  &1 &-i&-1&i
\end{array}
\end{equation}
\vspace{0.4cm}

{\bf Table 3.3}. Characters and multiplicities in the reduction $S(4)>C_4$.

\begin{equation}
\label{t20}
\begin{array}{r|rrrr|l}
k & e& (1,2,3,4) &(1,3)(2,4)& (4,3,2,1)&m(f,0)\\
\cline{1-6}
\chi^{[4]}(k) & 1& 1& 1&1&1\\ 
\chi^{[31]}(k) & 3& -1& -1&-1&0\\
\chi^{[22]}(k) & 2& 0& 0&2&1\\
\chi^{[211]}(k) & 3&1&1&-1&1\\
\chi^{[1111]}(k)& 1&-1&-1&1&0\\
\chi^0(k) & 1& 1& 1&1\\
\end{array}
\end{equation}

{\bf Table 3.4}. The multiplicities $m((l,\kappa),f)$ of representations $D^f$ of $S(4)$ and  $m((\lambda,\kappa),0)$ of $C_4$-periodic eigenmodes in the general reduction $O(3)>S(4)>C_4$ for $(l, \kappa), l \leq 4$. Because of eq. \ref{t22}, only the pairs $(l,\kappa)$ with $\kappa= (-1)^l$ 
contribute to the $C_4$-periodic modes. 
\begin{equation}
\label{t23a}
\begin{array}{l|lllll|l}
(l,\kappa),\:f:& \left[4\right]&\left[31\right]&\left[22\right]
&\left[211\right]&\left[1111\right]&m((l,\kappa),0)\\
\cline{1-7}
(0,1)     & 1             &0              &0
&0               &0       &1\\
(1,-1)    & 0             &1              &0
&0               &0       &0\\
(2,1)     & 0             &1              &1
&0        & 0             &1\\
(3,-1)    &1              &1              &0
&1        &0              &2\\
(4,1)     &1              &1              &1
&1               &0       &3\\

\end{array}
\end{equation}

\section{The tetrahedral 3-simplex ${\cal S}_0(3)$ on the sphere $S^3$.}

We consider the 3-sphere $S^3<E^4$ and an inscribed regular 4-simplex with
its vertices enumerated as $1,2,3,4,5$. The full point symmetry of the 4-simplex is 
$S(5)$. Central projection of the 3-faces of this simplex to $S^3$
yields a tiling with $5$ tetrahedral tiles. We choose the tetrahedron obtained by dropping 
the vertex $5$ as the simplicial manifold ${\cal S}_0(3)$. Its internal point symmetry group 
is $S(4)$.  The homotopy group $\pi_1({\cal S}_0(3))$ of the Platonic tetrahedron is described by Everitt in 
\cite{EV04} by a graph algorithm. Its prime dimension $5$ identifies it and the group of deck transformations as the cyclic group
$C_5$. In the next subsections we implement the group/subgroup analysis 
in analogy to the previous sections with the goal to characterize 
the harmonic analysis on ${\cal S}_0(3)$.

\subsection{The reduction $S(5)>C_5$.}
The characters of the irreducible representations 
of $S(5)$ for all $7$ classes and irreducible representations from \cite{HA62} p. 276 are given in {\bf Table 4.1}.

We shall need the characters of Table 4.1 for the reduction of the rotation group.

The cyclic group $C_5$ has the elements
\begin{equation}
\label{t32}
 g=(1,2,3,4,5), g^2=(1,3,5,2,4), g^3=(1,4,2,5,3), g^4=(1,,5,4,3,2), g^5=e. 
 \end{equation}
They belong to the classes $(4)(1)$ or $(5)$ of $S(5)$. The computation 
of the multiplicity $m(f,0)$ of the identity representation $D^0(C_5)$ is straightforward and we include it in the last column of {\bf Table 4.1}. The representation $D^0$ is contained once in the representations 
$f=\left[ 5\right], \left[ 11111\right],\left[ 32\right],\left[ 221\right]$,
twice in the representation $f=\left[ 311\right]$, but not in the representations 
$f=\left[ 41\right], f=\left[ 2111\right]$. Again we must determine the corresponding subspaces
for the first set of representations.

\subsection{$O(4,R)$ and Weyl reflections acting on $S^3$.}

We shall adopt from \cite{KR05} the coordinate description  of the sphere 
$S^3$, equivalent to an element $u \in SU(2,C)$, and relate it to the  vector notation 
$x=(x_0,x_1,x_2,x_3)$ in $E^4$ by 
\begin{eqnarray}
\label{t33}
&&u= \left[
\begin{array}{ll}
z_1& z_2\\
-\overline{z}_2& \overline{z}_1
\end{array}
\right],\: {\rm det}(u) = z_1\overline{z}_1+z_2\overline{z}_2=1,
\\ \nonumber 
&&z_1= x_0-ix_3,\: z_2=-(x_2+ix_1).
\end{eqnarray}
We use  the isomorphism
$SO(4,R) \sim (SU^l(2,C)\times SU^r(2,C))/Z_2$ with the action defined as 
\begin{eqnarray}
 \label{t34}
&&(SU^l(2,C) \times SU^r(2,C)) \times S^3 \rightarrow S^3, 
\\ \nonumber 
&&(g_l, g_c) \times u \rightarrow (g_l)^{-1} u g_r.
\end{eqnarray}
Here $Z_2$ denotes the subgroup generated by the action of 
$(-e,-e) \in (SU^l(2,C) \times SU^r(2,C))$ in eq. \ref{t34} which preserves 
any point $u \in S^3$. We stress that the half-angular parameters  here  must be
used in their full range appropriate for $SU(2)$.
The spherical harmonics on $S^3$ are shown in \cite{KR05} to be equal
to the Wigner representation functions $D^j$ of $SU(2)$ eq. \ref{a2a},
\begin{eqnarray}
 \label{t35}
&&Y^j_{m_1m_2}(u) = D^j_{m_1m_2}(u),\: j=0,1/2,1,3/2,..,\: -j\leq m_i\leq j,\: i=1,2,
\\ \nonumber 
&&\int \overline{D}^{j'}_{m_1'm_2'}(u)\: D^j_{m_1m_2}(u)\: d\mu(u)
= \delta_{j',j}\delta_{m_1',m_1}\delta_{m_2',m_2}.
\end{eqnarray}
The measure $d\mu(u)$ can be expressed by three Euler angle parameters \cite{ED57} pp. 62-3.
We denote the irreducible representation obtained by the action of rotations and
Weyl reflections on the spherical harmonics eq. \ref{t35} as $D^{(j,j)}$. 
In the appendix we summarize some properties of the Wigner $D^j$ functions.
In terms of the complex numbers $(z_1,z_2,\overline{z}_1,\overline{z}_2)$ in eq. \ref{t33}, the spherical harmonics are homogeneous polynomials of degree $(2j)$, see eq. \ref{a2a}. 
 Under the action  
of rotations they transform, due to their  representations properties, irreducibly with 
a representation $D^{(j,j)}$ given by direct product matrices, 
\begin{eqnarray}
 \label{t36}
&&(T_{(g_l,g_r)}D^j_{m_1m_2})(u)
=D^j_{m_1,m_2}(g_l^{-1}ug_r)
\\ \nonumber 
&&=\sum_{m_1',m_2'} D^j_{m_1'm_2'}(u)\left[D^j_{m_1m_1'}(g_l^{-1}) D^j_{m_2'm_2}(g_r)\right],
\\ \nonumber
&&D^{(j,j)}_{m_1'm_2',m_1m_2}(g_l,g_r)
= \langle (j,j)(m_1'm_2')|T_{(g_l,g_r)}|(j,j)(m_1m_2)\rangle
\\ \nonumber
&&=D^j_{m_1m_1'}(g_l^{-1}) D^j_{m_2'm_2}(g_r) 
\end{eqnarray}
where finally we passed for convenience to a bracket notation.
The defining action of $SO(4,R)$ on $E^4$ by real $4\times 4$ rotation matrices is included in this notation as the representation $D^{(\frac{1}{2},\frac{1}{2})}$.

For  the character of the irreducible representation $D^{(j,j)}$ of a rotation
it follows from eq. \ref{t36} that it 
is given  by
\begin{equation}
\label{t37}
\chi^{(j,j)}(T_{(g_l,g_r)}) = {\rm Trace}(D^j(g_l^{-1}) \times D^j(g_r))= \chi^j(g_l^{-1})\chi^j(g_r), 
\end{equation}
in terms of pairs of  characters $\chi^j(g)$ of $SU(2,C)$. The dimension 
of $D^{(j,j)}$ from eq. \ref{t37} is $(\chi^j(e))^2=(2j+1)^2$.

To determine the action of the symmetric group $S(5)$ on $S^3$ we need the action of a Weyl reflection eq. \ref{t3a} in the coordinate description eq. \ref{t33} and at the same time must extend the action from
$SO(4,R)$ to $O(4,R)$. We start with the 
Weyl reflection 
for the particular Weyl vector $a_0=(1,0,0,0)$ which acts as 
\begin{equation}
\label{t38}
(W_{a_0}\times  (x_0,x_1,x_2,x_3)) \rightarrow (-x_0,x_1,x_2,x_3).
\end{equation} 
Inserting this reflection into the complex coordinates of $u$ from equation \ref{t33} we find that 
the image of $u$ under $W_{a_0}$ can be rewritten in matrix notation as
\begin{equation}
 \label{t39}
W_{a_0}(u) = -u^{\dagger}.
\end{equation}
In addition we need the representations and characters for reflections
from $O(4,R)$. 

For $u \in SU(2,C)$ we need the properties  eqs. \ref{a1} and \ref{a4}.
This equation and the homogeneity eq. \ref{a5} of degree $2j$ allows us to rewrite the Weyl operator 
$T_{a_0}$ acting on a spherical harmonic eq.\ref{t35} as
\begin{eqnarray}
 \label{t41}
&&(T_{a_0} D^j_{m_1m_2})(u)= D^j_{m_1m_2}(-u^{\dagger})
\\ \nonumber
&&= (-1)^{2j}D^j_{m_2m_1}(q^{-1}uq)= (-1)^{2j}(T_{(q,q)}D^j_{m_2m_1})(u).
\end{eqnarray}
We used the homogeneity eq. \ref{a5} to extract the phase factor $(-1)^{2j}$.
With respect to an action eq. \ref{t36} of $SO(4,R)$ following the reflection $W_{a_0}(u)$, we get from eq. \ref{t41} 
\begin{equation}
 \label{t42}
((T_{(g_l,g_r)} T_{a_0})D^j_{m_1,m_2})(u) 
= (-1)^{2j} ((T_{(g_l,g_r)} T_{(q,q)})D^j_{m_2,m_1})(u)
\end{equation}
For the Weyl operator $T_{a_0}$ eq. \ref{t41} we find the involutive and conjugation  properties
\begin{eqnarray}
 \label{t43}
&&T_{a_0} T_{a_0}= I,
\\ \nonumber 
&&T_{a_0} T_{(g_l,g_r)} T_{a_0}= T_{(g_r,g_l)} 
\end{eqnarray}

The Weyl vector $a_0$ in the coordinates eq. \ref{t33} corresponds to $u=e$.
For a general unit Weyl vector $a$ we can choose  a rotation $R \in SO(4,R)$ such that
$a= Ra_0$. Then by the general operator relation between Weyl reflections and rotations 
\begin{equation}
\label{t44}
W_a= W_{Ra_0} =  R W_{a_0}R^{-1},
\end{equation}
the general Weyl reflection $W_a$ can be expressed in terms of the particular reflection $W_{a_0}$. It remains to construct a rotation $R$. In terms of the coordinates eq. 
\ref{t33} we note that
\begin{equation}
\label{t45}
((e, v) \times e)=v
\end{equation}
which means that the element $(e, v)\in (SU^l(2,C) \times SU^r(2,C))$ applied to
$a_0$ gives the most general point $ a\in S^3$  and so provides the rotation $R$. 

With this choice of $R$, the action of a general Weyl reflection operator on a spherical harmonic  can be written by use of eq. \ref{t43} as 
\begin{equation}
 \label{t46}
T_{a}= T_{(e,v^{-1}_a)} T_{a_0} T_{(e,v_a)} = T_{(v_a,v^{-1}_a)} T_{a_0}
\end{equation}
Here the matrix $v_a$ is obtained from  the Weyl vector $a$ by inserting 
the vector components of $x=a$ into the matrix eq. \ref{t33}.
From the factorization eq. \ref{t46}  of the general Weyl reflection operator 
and from eq. \ref{t36} we find the matrix representation of $T_a$ as
\begin{eqnarray}
 \label{t47}
&&\langle (j,j)(m_1'm_2')| T_{(v_a,v^{-1}_a)}T_{a_0}|(j,j) (m_1m_2) \rangle 
\\ \nonumber 
&&=(-1)^{2j}
\langle (j,j)(m_1'm_2')| T_{(v_a,v^{-1}_a)} T_{(q,q)}|(j,j) (m_2m_1) \rangle 
\\ \nonumber
&&=(-1)^{2j} D^j_{m_2m_1'}(q^{-1}v^{-1}_a)D^j_{m_2'm_1}(v^{-1}_aq).
\end{eqnarray}
Here we used the complex conjugation property $\overline{u}=q^{-1}uq$ of  $SU(2,C)$, eq. \ref{a1}. 
The trace of this representation becomes 
\begin{equation}
\label{t48}
\chi^{(j,j)} (T_a)= (-1)^{2j} \chi^j(q^{-1}v^{-1}_a(-q)(v^{-1}_a)^T)
=\chi^j(\overline{v^{-1}_a} (v^{-1}_a)T) =\chi^j(e)=2j+1.
\end{equation}
This result is independent of $v$ as  expected from eq. \ref{t46} since 
all Weyl operators by eq. \ref{t44} are conjugate to $T_{a_0}$ and the trace is independent
of conjugations.

All operators which are not pure rotations can be written with the help of 
eq. \ref{t43} in the form
\begin{equation}
 \label{t49}
T= T_{(g_l,g_r)} T_{a_0}.
\end{equation}
The matrix elements of this operator are given similar to eq. \ref{t47}
by
\begin{eqnarray}
\label{t49a}
&&(-1)^{2j} 
\langle (j,j) (m_1'm_2') | T_{(g_lq,g_rq)}|(j,j) (m_2m_1)\rangle
\\ \nonumber
&&=(-1)^{2j} D^j_{m_2m_1'}(q^{-1}g_l^{-1})D^j_{m_2'm_1}(g_rq).
\end{eqnarray}
with $q$ given in eq. \ref{a1}. It follows from eq. \ref{t49a} that the character of this operator in the representation 
$D^{(j,j)}$ is given by 
\begin{eqnarray}
\label{t50}
&&\chi^{(j,j)}(T_{(g_l,g_r)} T_{a_0})= (-1)^{2j}  
\sum_{m_1m_2} D^j_{m_2m_1}(q^{-1}g_l^{-1}) D^j_{m_2m_1}(g_rq)
\\ \nonumber
&&=\chi^j(\overline{g_l^{\dagger}}g_r^T)=\chi^j((g_rg_l)^T)= \chi^j(g_rg_l) 
\end{eqnarray}
Here we used from eq. \ref{a1} $q^T=(-q)$, which leads to the cancellation of the 
phase factor $(-1)^{2j}$, and $q^{-1}g_l^{-1}q= q^{-1}g_l^{\dagger}q=g_l^T$.

\subsection{The representations of $S(5)$.}

From the Coxeter diagram of $S(5)$ we can construct in $E^4$ a set of $4$ Weyl vectors 
associated with  the $4$ generators. The scalar products must form, 
\cite{HU90} pp. 108-10,
the matrix
\begin{equation}
\label{t51} 
 \langle a_i, a_j \rangle 
= \left[
\begin{array}{llll}
1& \frac{1}{2}& 0&0\\
\frac{1}{2}&1& \frac{1}{2}&0\\
0&\frac{1}{2}&1&\frac{1}{2}\\
0&0&\frac{1}{2}&1
\end{array}
\right]
\end{equation}
Here the off-diagonal scalar products  are determined from the Coxeter diagram 
by $\langle a_i,a_{i+1}\rangle = \cos(\pi/3)=\frac{1}{2},\; i=1,2,3$.
We fulfill eq. \ref{t51} by the choice   of the four Weyl vectors given in {\bf Table 4.2}.

From eq. \ref{t46} we compute for each Weyl vector $w_i, i=1,\ldots,4$ 
the  matrices $v_i$ which appear in the reflection operator and include them in 
the last column of {\bf Table 4.2}. 

From the generators of $S(5)$ we can construct representatives $g_k$ of the
seven classes $k$ appearing in {\bf Table 4.2} and express them as products of
the Weyl reflection operators $W_i:=W_{a_i}$ in {\bf Table 4.3}.

For products of up to $4$ Weyl reflection operators eq. \ref{t46} we find
with eq. \ref{t43} 
\begin{eqnarray}
 \label{t56}
&&T_{a}= T_{(v_a,v_a^{-1})} T_{a_0},
\\ \nonumber
&&T_{b} T_{a} = T_{(v_bv_a^{-1},v_b^{-1} v_a)},
\\ \nonumber
&&T_{c}T_{b}T_{a}=T_{(v_cv_b^{-1}v_a,v_c^{-1} v_bv_a^{-1})} T_{a_0},
\\ \nonumber 
&&T_{d}T_{c}T_{b}T_{a}
=T_{(v_dv_c^{-1}v_bv_a^{-1},v_d^{-1} v_cv_b^{-1} v_a)}.
\end{eqnarray}
Here products with an even number of Weyl operators are rotations, those  with an odd number become a product of
a rotation with $T_{a_0}$.

\subsection{Multiplicity in the reduction $O(4,R)>S(5)>C_5$.}
To write  explicitly  the action of the class representatives from {\bf Table 4.3} 
we must convert the corresponding   product of Weyl operators for each class representative with the help
of eq. \ref{t56} into an explicit action of $SU^l(2,C) \times SU^r(2,C)$
on $u$ or on $(-u^{\dagger})$ respectively. Once we have computed 
these actions, 
the characters in the representation $ D^{(j,j)}$  are found  
by use of eqs. \ref{t36}, \ref{t50} respectively. 
It turns out that for five classes $k\in S(5)$ the characters 
$\chi^{(j,j)}(k)$ are periodic with periods $2,3,4,5$ respectively 
wrt. $(2j)$.

The results and corresponding recursion relations are given in  {\bf Tables 4.7, 4.8}.

The real scalar product 
\begin{equation}
 \label{t57}
m((j,j),f) = \frac{1}{5!}\sum_{k \in S(5)}  n(k)\:\chi^{(j,j)}(k) \chi^f(k),
\end{equation}
with $n(k)$ the number of elements in class  $k \in S(5)$, then yields the multiplicity $m((j,j),f)$ in the reduction $O(4,R)>S(5)$.
In {\bf Table 4.9} we write down the multiplicities computed according to eq. \ref{t57} 
for the reduction of representations in $O(4,R)>S(5)$. The multiplicity 
of $C_5$-periodic modes of $O(4,R)$ is then given by
\begin{equation}
\label{t57a}
m((j,j),0)=\sum_f m((j,j),f)\:m(f,0).
 \end{equation}
in the last column of {\bf Table 4.9}.

Of the $506$ harmonic polynomials
for $O(4,R)$ up to degree $2j=10$, $101$ are $C_5$-periodic eigenmodes of ${\cal S}_0(3)$.
and so demonstrate on average
the topological selection rules for the simplicial manifold ${\cal S}^0(3)$.
There are other specific selection rules like the absence of $C_5$-periodic eigenmodes for
$2j=1$.
By use of eq. \ref{t57} and the characters given in {\bf Tables  4.7, 4.8}, {\bf Table 4.9} can be 
extended to any value $(2j) > 10$.
We note another recursion relation for the multiplicity $m((j,j),f)$: for all except
the first two classes, the characters $\chi^{(j,j)}(k)$ from {\bf Tables 4.7, 4.8} obey
\begin{equation}
\label{t57d}
k=(3)(1)^2, (2)^2(1),(3)(2),(4)(1),(5):\: \chi^{(j+30,j+30)}(k)=\chi^{(j,j)}(k). 
\end{equation}
Using in eq. \ref{t57} the characters for the first two classes
from {\bf Table 4.9} and for the others the periodicity  eq. \ref{t57d} one finds for the multiplicity 
the recursion relation 
\begin{equation}
\label{t57c}
m((j+30,j+30),f)= m((j,j),f)+(2j+36). 
\end{equation}
\vspace{0.3cm}

\subsection{The explicit reduction $S(5)>C_5$.}
For the partitions $f=\left[32\right],\left[221\right],\left[311\right]$ 
we now compute the $C_5$-periodic states $\psi^{f 0}$ as linear combinations of the basis states $\phi_j$ in the Young orthogonal
representation. For each representation  $D^f$ which admits $D^0$, they are given as the single or at most two eigenvectors of eigenvalue $1$ 
for the Coxeter element generating  $C_5$ in the representations spaces for those partitions which from {\bf Table 4.1} reduce to $D^0$.
We give as a vector the coefficients $\psi^{f 0}_j$ in the 
linear combinations
\begin{equation}
\label{t58a}
\psi^{f 0}= \sum_j \phi_j\psi^{f 0}_j. 
\end{equation}
The representation matrices for the generators of $S(5)$ and $C_5$ and the coefficients
in eq. \ref{t58a} are given 
in {\bf Tables 4.4-4.6}.

\subsection{Tables for $O(4,R)>S(5)>C_5$.}

{\bf Table 4.1}. Characters for the group $S(5)$ for irreducible representations $D^f$ with 
partition $f$, 
classes $k$ in cycle notation, number $n(k)$ 
of elements in class $k$, and multiplicity $m(f,0)$ of the
identity representation $D^0(C_5)$.

\begin{equation}
\label{t31}
\begin{array}{r|rrrrrrr|l}
\chi^f(k) & (1)^5&(2)(1)^3&(2)^2(1)&(3)(1)^2&(3)(2)&(4)(1)&(5)&m(f,0)\\
\cline{1-9}
&&&&&&&&\\
n(k) &1&10&15&20&20&30&24\\
\left[5\right]       &1&1&1&1&1&1&1&1\\
\left[11111\right]   &1&-1&1&1&-1&-1&1&1\\
\left[41\right]      &4&2&0&1&-1&0&-1&0\\
\left[2111\right]    &4&-2&0&1&1&0&-1&0\\
\left[32\right]      &5&1&1&-1&1&-1&0&1\\
\left[221\right]     &5&-1&1&-1&-1&1&0&1\\
\left[311\right]     &6&0&-2&0&0&0&1&2\\
\end{array}.
 \end{equation}

{\bf Table 4.2}. Weyl vectors $a_i$ and matrices $v_i$ in the Weyl operators
eq. \ref{t46} for the generators $(i,i+1)$ of the group $S(5)$.

\begin{equation}
\label{t52}
\begin{array}{l|l|l|l}
i&(i,i+1)& a_i& v_i\\
\hline
&&&\\
1&(1,2)& (0,0,0,1)&
\left[
\begin{array}{ll}
-i&0\\
0&i
\end{array}
\right]\\
&&&\\
2&(2,3)& (0,0,\sqrt{\frac{3}{4}}, \frac{1}{2})&
\left[
\begin{array}{ll}
-\frac{i}{2}&-\frac{1}{2}\sqrt{3}\\
\frac{1}{2}\sqrt{3}&\frac{i}{2}
\end{array}
\right]\\
&&&\\
3&(3,4)& (0, \sqrt{\frac{2}{3}},\sqrt{\frac{1}{3}},0)&
\left[
\begin{array}{ll}
0&-(\sqrt{\frac{1}{3}}+i\sqrt{\frac{2}{3}})\\
(\sqrt{\frac{1}{3}}-i\sqrt{\frac{2}{3}})& 0
\end{array}
\right]\\
&&&\\
4&(4,5)& (\sqrt{\frac{5}{8}},\sqrt{\frac{3}{8}}, 0,0)&
\left[
\begin{array}{ll}
\sqrt{\frac{5}{8}}&-i\sqrt{\frac{3}{8}}\\
-i\sqrt{\frac{3}{8}}&\sqrt{\frac{5}{8}}
\end{array}
\right]\\
\end{array}
\end{equation}

{\bf Table 4.3}. Representatives for the seven classes $k$ of $S(5)$ in terms 
of Weyl reflections $W_i:=W_{a_i}$, and matrices $g_l,g_r, g_rg_l$ appearing in 
eqs. \ref{t37}, \ref{t49}, \ref{t50}.

\begin{equation}
\label{t53}
\begin{array}{l|llllll}
k: &(1)^5&(2)(1)^3&(2)^2(1)  &(3)(1)^2&&\\
\hline
g_k&e    &(1,2)   &(1,2)(3,4)&(1,2)(2,3)&&\\
&I& W_1& W_1W_3& W_1W_2&&\\
g_l
&
&
&
\left[
\begin{array}{ll}
0&\frac{\sqrt{2}-i}{\sqrt{3}}\\
-\frac{\sqrt{2}+i}{\sqrt{3}}&0
\end{array}
\right]
&
\left[
\begin{array}{ll}
\frac{1}{2}&-\frac{i\sqrt{3}}{2}\\
 -\frac{i\sqrt{3}}{2}&\frac{1}{2}
\end{array}
\right]
\\
g_r
&
&
&
\left[
\begin{array}{ll}
0&\frac{\sqrt{2}-i}{\sqrt{3}}\\
-\frac{\sqrt{2}+i}{\sqrt{3}}&0
\end{array}
\right]
&
\left[
\begin{array}{ll}
\frac{1}{2}&-\frac{i\sqrt{3}}{2}\\
 -\frac{i\sqrt{3}}{2}&\frac{1}{2}
\end{array}
\right]
\\
g_rg_l
&
&
&
&
\\
\end{array}
\end{equation}
\begin{equation}
\begin{array}{l|llll}
k: &(3)(2)         &(4)(1)         &(5)                 &\\
\hline
g_k&(1,2)(2,3)(4,5)&(1,2)(2,3)(3,4)&(1,2)(2,3)(3,4)(4,5)&\\
   &W_1W_2W_4      & W_1W_2W_3     &W_1W_2W_3W_4&\\
g_l
&
&
&
\left[
\begin{array}{ll}
\frac{2-2\sqrt{5}-i(\sqrt{2}+\sqrt{10})}{8}
&\frac{3\sqrt{2}-\sqrt{10}+i(-6-2\sqrt{5})}{8\sqrt{3}}\\
\frac{-3\sqrt{2}+\sqrt{10}+i(-6-2\sqrt{5})}{8\sqrt{3}}
&\frac{2-2\sqrt{5} +i(\sqrt{2}+\sqrt{10})}{8} 
\end{array}
\right]
&
\\
g_r
&
&
&
\left[
\begin{array}{ll}
\frac{2+2\sqrt{5}+i(-\sqrt{2}+\sqrt{10})}{8}
&\frac{3\sqrt{2}+\sqrt{10}+i(-6+2\sqrt{5})}{8\sqrt{3}}\\
\frac{-3\sqrt{2}-\sqrt{10}+i(-6+2\sqrt{5})}{8\sqrt{3}}
&\frac{2+2\sqrt{5} +i(\sqrt{2}-\sqrt{10})}{8} 
\end{array}
\right]
&
\\
g_rg_l
&
\left[
\begin{array}{ll}
-\frac{1}{2}&-\frac{i\sqrt{3}}{2}\\
 -\frac{i\sqrt{3}}{2}&-\frac{1}{2}
\end{array}
\right]
&
\left[
\begin{array}{ll}
-\frac{i}{\sqrt{2}}&-\frac{\sqrt{2}+2i}{2\sqrt{3}}\\
 \frac{\sqrt{2}-2i}{2\sqrt{3}}&\frac{i}{\sqrt{2}}
\end{array}
\right]
&
&
\\
\end{array}
\end{equation}

{\bf Table 4.4}. Periodic state for the partition $f=\left[32\right]$.
Young tableau basis:
\begin{eqnarray}
\label{t58}
&&1=
\left[ 
\begin{array}{lll}
1&2&3\\ 
4&5&
\end{array}
\right],
2=
\left[ 
\begin{array}{lll}
1&3&4\\ 
2&5&
\end{array}
\right],
3=
\left[ 
\begin{array}{lll}
1&2&4\\ 
3&5&
\end{array}
\right],
\\ \nonumber
&&
4=
\left[ 
\begin{array}{lll}
1&3&5\\ 
2&4&
\end{array}
\right],
5=
\left[ 
\begin{array}{lll}
1&2&5\\ 
3&4&
\end{array}
\right].
\end{eqnarray}
Representation $D^{\left[ 32\right]}$ for generators from
\cite{HA62} pp. 227-8:
\begin{eqnarray}
\label{t59}
&& D^{[32]}(1,2)= 
\left[
\begin{array}{lllll}
1&0&0&0&0\\
0&-1&0&0&0\\
0&0&1&0&0\\
0&0&0&-1&0\\
0&0&0&0&1 
\end{array}
\right],\:
D^{[32]}(2,3)= 
\left[
\begin{array}{lllll}
1&0&0&0&0\\
0&\frac{1}{2}&\sqrt{\frac{3}{4}}&0&0\\
0&\sqrt{\frac{3}{4}}&-\frac{1}{2}&0&0\\
0& 0&0&\frac{1}{2}&\sqrt{\frac{3}{4}}\\
0& 0&0&\sqrt{\frac{3}{4}}&-\frac{1}{2}
\end{array}
\right],
\\ \nonumber
&& D^{[32]}(3,4))= 
\left[
\begin{array}{lllll}
-\frac{1}{3}&0&\sqrt{\frac{8}{9}}&0&0\\
0&1&0&0&0\\
\sqrt{\frac{8}{9}}&0&\frac{1}{3}&0&0\\
0&0&0&-1&0\\
0&0&0&0&1 
\end{array}
\right],\:
D^{[32]}(4,5)= 
\left[
\begin{array}{lllll}
1&0&0&0&0\\
0&-\frac{1}{2}&0&\sqrt{\frac{3}{4}}&0\\
0&0&-\frac{1}{2}&0&\sqrt{\frac{3}{4}}\\
0&\sqrt{\frac{3}{4}}&0&\frac{1}{2}&0\\
0&0&\sqrt{\frac{3}{4}}&0&\frac{1}{2} 
\end{array}
\right].
\end{eqnarray}
The matrix  $D^{\left[32\right]}(4,5)$ given in \cite{HA62} p. 228 has an error  and was replaced.

Representation of Coxeter element of $C_5$:
\begin{equation}
\label{t60}
D^{\left[32\right]}((1,2)(2,3)(3,4)(4,5))=
\left[
\begin{array}{lllll}
-\frac{1}{3}&0&-\sqrt{\frac{2}{9}}&0&\sqrt{\frac{2}{3}}\\ 
-\sqrt{\frac{2}{3}}&\frac{1}{4}&\sqrt{\frac{1}{48}}&-\sqrt{\frac{3}{16}}&-\frac{1}{4}\\
-\sqrt{\frac{2}{9}}&-\sqrt{\frac{3}{16}}&\frac{1}{12}&\frac{3}{4}&-\sqrt{\frac{1}{48}}\\
0&\sqrt{\frac{3}{16}}&-\frac{3}{4}&\frac{1}{4}&-\sqrt{\frac{3}{16}}\\
0&-\frac{3}{4}&-\sqrt{\frac{3}{16}}&-\sqrt{\frac{3}{16}}&-\frac{1}{4}
\end{array}
\right]
\end{equation}
Eigenvector of Coxeter element  for the single eigenvalue $1$ eq. \ref{t58a}:
\begin{equation}
\label{t61}
\psi^{\left[32\right] 0}=(\sqrt{\frac{2}{3}},-1,-\sqrt{\frac{1}{3}},-\sqrt{\frac{1}{3}},1). 
\end{equation}

{\bf Table 4.5}. Periodic state for the partition $f=\left[221\right]$.
The Young tableaus of the basis are the mirror images of those for the associate 
partition $\left[32\right]$ eq. \ref{t58}. The representation matrices can 
then be computed from those of $D^{\left[ 32\right]}$ eq. \ref{t59}.

Representation of the Coxeter element:
\begin{equation}
\label{t62}
D^{\left[221\right]}((1,2)(2,3)(3,4)(4,5))=
\left[
\begin{array}{lllll}
-\frac{1}{3}&0&\sqrt{\frac{2}{9}}&0&\sqrt{\frac{2}{3}}\\ 
-\sqrt{\frac{2}{3}}&\frac{1}{4}&-\sqrt{\frac{1}{48}}&\sqrt{\frac{3}{16}}&-\frac{1}{4}\\
\sqrt{\frac{2}{9}}&\sqrt{\frac{3}{16}}&\frac{1}{12}&\frac{3}{4}&\sqrt{\frac{1}{48}}\\
0&-\sqrt{\frac{3}{16}}&-\frac{3}{4}&\frac{1}{4}&\sqrt{\frac{3}{16}}\\
0&-\frac{3}{4}&\sqrt{\frac{3}{16}}&\sqrt{\frac{3}{16}}&-\frac{1}{4}
\end{array}
\right]
\end{equation}
Eigenvector of the Coxeter element for the eigenvalue $1$:
\begin{equation}
\label{t63}
\psi^{\left[221\right] 0}=(\sqrt{\frac{2}{3}},-1,\sqrt{\frac{1}{3}},\sqrt{\frac{1}{3}},1). 
\end{equation}

{\bf Table 4.6}. Periodic state for the partition $f=\left[311\right]$ in eq. \ref{t58a}:

Young tableau basis: 
\begin{eqnarray}
\label{t64}
&&1=\left[\begin{array}{lll}
1&2&3\\
4&&\\
5&&
\end{array}
\right],\:
2=\left[\begin{array}{lll}
1&2&4\\
3&&\\
5&&
\end{array}
\right],\:
3=\left[\begin{array}{lll}
1&3&4\\
2&&\\
5&&
\end{array}
\right],\:
\\ \nonumber
&&4=\left[\begin{array}{lll}
1&2&5\\
3&&\\
4&&
\end{array}
\right],\:
5=\left[\begin{array}{lll}
1&3&5\\
2&&\\
4&&
\end{array}
\right],\:
6=\left[\begin{array}{lll}
1&4&5\\
2&&\\
3&&
\end{array}
\right],\:
\end{eqnarray}
The representation matrices of the generators are given in \cite{HA62}
pp. 229-30. Computed  matrix representation of the Coxeter element: 
\begin{equation}
\label{t65}
D^{\left[311\right]}((1,2)(2,3)(3,4)(4,5))=
\left[
\begin{array}{llllll}
\frac{1}{3}&-\sqrt{\frac{1}{18}}&0
&\sqrt{\frac{5}{6}}&0&0\\
\sqrt{\frac{2}{9}}&\frac{1}{24}&-\sqrt{\frac{3}{64}}
&-\sqrt{\frac{5}{192}}&\sqrt{\frac{45}{64}}&0\\
\sqrt{\frac{2}{3}}&\sqrt{\frac{1}{192}}&\frac{1}{8}
&-\sqrt{\frac{5}{64}}&-\sqrt{\frac{15}{64}}&0\\
0&\sqrt{\frac{15}{64}}&-\sqrt{\frac{5}{64}}
&\frac{1}{8}&-\sqrt{\frac{1}{192}}&\sqrt{\frac{2}{3}}\\
0&\sqrt{\frac{45}{64}}&\sqrt{\frac{5}{192}}
&\sqrt{\frac{3}{64}}&\frac{1}{24}&-\sqrt{\frac{2}{9}}\\
0&0&\sqrt{\frac{5}{6}}
&0&\sqrt{\frac{1}{18}}&\frac{1}{3}
\end{array}
\right]
\end{equation}
Two eigenvectors in eq. \ref{t58a} for eigenvalue $1$ of eq. \ref{t65}:
\begin{eqnarray}
\label{t66}
&&q^{\left[311\right] 0,1}
=(\sqrt{\frac{49}{45}},\sqrt{\frac{2}{45}},\sqrt{\frac{8}{15}},
\sqrt{\frac{2}{3}},0,1),
\\ \nonumber
&&q^{\left[311\right] 0,2}
=(\sqrt{\frac{8}{45}},\sqrt{\frac{49}{45}},-\sqrt{\frac{1}{15}},
\sqrt{\frac{1}{3}},1,0),
\\ \nonumber
\end{eqnarray}

{\bf Table 4.7}. The angles $\phi$ and characters $\chi^{(j,j)}(k)$ for the classes $k \in S(5)$ 
and the representations $D^{(j,j)}$ of $O(4,R)$ with $(2j)=0,1,2,3,4,5$. General 
expressions and recursion relations are given in {\bf Table 4.8}.

\begin{equation}
\label{t56a} 
\begin{array}{l|l|lllllll}
k&\phi&\chi^{(j,j)}(k),(2j):& 0&1&2&3&4&5 \\
\hline
(1)^5&\phi(g_l)/2=\phi(g_r)/2=0&&1&4&9&16&25&36\\
(2)(1)^3&\phi(g_rg_l)/2=0&&1&2&3&4&5&6\\
(3)(1)^2&\phi(g_l)/2=\phi(g_r)/2=\pi/3&&1&1&0&1&1&0\\
(2)^2(1)&\phi(g_l)/2=\phi(g_r)/2=\pi/2&&1&0&1&0&1&0\\
(3)(2)&\phi(g_rg_l)/2=2\pi/3&&1&-1&0&1&-1&0\\
(4)(1)&\phi(g_rg_l)/2=\pi/2&&1&0&-1&0&1&0\\
(5)&\phi(g_l)/2=2\pi/10,\phi(g_r)/2=6\pi/10&&1&-1&-1&1&0&1\\
\end{array}
\end{equation}

{\bf Table 4.8} General expressions and recursion relations for the
characters $\chi^{(j,j)}(k)$ as functions of classes $k$ of $S(5)$. The initial values are given in  {\bf Table 4.7}.

\begin{equation}
 \label{t56b}
\begin{array}{l|l}
k&\chi^{(j,j)}(k)\\
\hline
(1)^5&\chi^{(j,j)}=(2j+1)^2\\
(2)(1)^3&\chi^{(j,j)}=(2j+1)\\
(3)(1)^2&\chi^{(j+3/2,j+3/2)}=\chi^{(j,j)}\\
(2)^2(1)&\chi^{(j+1,j+1)}=\chi^{(j,j)}\\
(3)(2)&\chi^{(j+3/2,j+3/2)}=\chi^{(j,j)}\\
(4)(1)&\chi^{(j+2,j+2)}=\chi^{(j,j)}\\
(5)&\chi^{(j+5/2,j+5/2)}=\chi^{(j,j)}
\end{array}
\end{equation}

{\bf Table 4.9}. Multiplicities $m((j,j),f)$ in the reduction of representations 
$D^{(j,j)}=\sum_j m((j,j),f)D^f$ from $O(4,R)$ to $S(5)$ as function of $(2j)=0,\ldots, 10$
and of all partitions $f$. The numbers
$m((j,j),0)$ in the last column denotes the total number of $C_5$-periodic modes 
for fixed $(2j)$, $\nu_0(f)$ in the last row the total number for a fixed partition $f$ up to $(2j)=10$.
\begin{equation}
\label{t70}
\begin{array}{r|rrrrrrrr|rr}
f&& \left[5\right]&\left[11111\right]&\left[41\right]&\left[2111\right]&\left[32\right]
&\left[221\right]&\left[311\right]&m((j,j),0)\\
\hline
&&&&&&&&\\
 (2j)   &   &  &  &  &  &  &  &\\
   
0   &&1 &  &  &  &  &  &  &1\\
1   &&  &  &1 &  &  &  &  &0\\
2   &&  &  &1 &  &1 &  &  &1\\
3   &&1 &  &1 &  &1 &  &1 &4\\
4   &&1 &  &2 &  &1 &1 &1 &5\\
5   &&1 &  &2 &  &2 &1 &2 &8\\
6   &&1 &  &3 &1 &3 &1 &2 &9\\
7   &&1 &  &4 &1 &3 &2 &3 &12\\
8   &&2 &  &4 &1 &4 &3 &4 &17\\
9   &&2 &  &5 &2 &5 &3 &5 &20\\
10  &&2 &1 &6 &2 &6 &3 &6 &24\\
\hline
&&&&&&&&\\
\nu_0(f)&&12&1&0&0&26&14&48&101 
\end{array}
\end{equation}
\vspace{0.3cm}
\subsection{Harmonic analysis on ${\cal S}_0(3)$.}
We summarize the basis construction for the harmonic analysis 
on ${\cal S}_0(3)$ in terms of $C_5$-periodic states on the sphere $S^3$.\\
({\bf i}) The spherical harmonics for fixed degree $2j=0,1,2,\ldots$ are the 
Wigner $D^j(u)$-functions given in eq. \ref{a2a}.\\ 
({\bf ii}) For the reduction $O(4,R)>S(5)$, the multiplicity of representations 
$D^f$ is known from {\bf Table 4.9} and computable from eq. \ref{t57} with 
the characters $\chi^{(j,j)}$ for $O(4,R)$ from {\bf Tables 4.7, 4.8} and $\chi^f$ for $S(5)$ 
given in {\bf Table 4.1}.  
The partitions $f=\left[ 41\right], \left[ 2111\right]$ are forbidden.
For allowed partitions $f$, the explicit state construction can employ the general Young operators from representation theory, \cite{WI59} pp. 112-23,
with the properties  
\begin{eqnarray}
\label{t71}
&&c^f_{r,s}:= \frac{|f|}{5!} \sum_{p'\in S(5)}
D^f_{r,s}(p') T_{p'},  
\\ \nonumber 
&&T_p c^f_{r,s}= \sum_{s'} c^f_{r',s} D^f_{r',r}(p)
\\ \nonumber
&& c^{\tilde{f}}_{\tilde{s},\tilde{r}} c^{f}_{r,s}
= \delta(\tilde{f},f)\delta(\tilde{r},r) c^{f}_{\tilde{s},s},
\\ \nonumber
&&(c^f_{r,s})^{\dagger}=c^f_{s,r}
\end{eqnarray}
The rows and columns $(r,s)$ of $D^f$ can be labelled by Young diagrams or by 
Yamanouchi symbols as in \cite{HA62}, compare for example eq. \ref{t27a}.
The second line in eq. \ref{t71}  shows that the application of a Young operator 
to a state yields a basis function for the irreducible representation $D^f$ 
with row lable $r$. To apply the Young operator to a spherical harmonic,
one needs the matrix elements of all permutations $p\in S(5)$ for the representation 
$D^{(j,j)}$ of $O(4,R)$. Any element $p \in S(5)$ can be written as a product of the
generators given in {\bf Table 4.2}. The corresponding product of operators $T_p$ is carried out with the help of 
eq. \ref{t56}. Then the expressions eq. \ref{t36} or \ref{t47} yield
the matrix elements $D^{(j,j)}(p)$. The matrix elements in the Young orthogonal
representation $D^f(p)$ 
are computed in a similar way from those of the generators, given 
for example in \cite{HA62} pp. 226-30, or constructed directly from 
\cite{HA62} pp. 219-20.\\
({\bf iii}) The reduction $S(5)>C_5$ for the 1-dimensional representations 
$f=\left[ 5\right],\left[ 11111\right]$ is unique, and for the 
partitions $f=\left[ 32\right],\left[ 221\right],\left[ 311\right]$
by eq. \ref{t58a} yields one or two unique states with the coefficients eqs. \ref{t61}, \ref{t63},
\ref{t66} respectively.\\
The $C_5$-periodic basis states on $S^3$ each belong to a fixed degree $(2j)$ from $O(4,R)$, carry an allowed 
partition $f$ from $S(5)$, are orthogonal with respect to $(2j)$ and $f$, and  stable under $C_5$.

\section{Conclusion.}
Methods of group theory allow to construct  and analyze the harmonic analysis 
on topological manifolds. This is demonstrated in section 4 for the simplicial manifold 
${\cal S}_0(3)$. 
The multiplicities provide 
the specific selection rules for the chosen simplex topology.
The symmetric group $S(5)$ plays a key role.
Its representations $f=\left[41\right],\left[2111\right]$
are eliminated from  the harmonic  analysis. 
The details for the basis construction are given in subsection 4.7.

In general, the harmonic analysis on two different 
manifolds  ${\cal M}, {\cal M}'$  covered by the sphere $S^ {(n-1)}$ 
is unified by the spherical harmonics and corresponding representations.
The differences between topologies 
appear in the form of different subgroups of deck transformations. In the harmonic analysis 
these  involve different  
group/subgroup representations and reductions in
$O(n,R)>deck({\cal M}),\: O(n,R)>deck({\cal M}')$.
Intermediate subgroups as $S(n+1)$ in eq. \ref{t1} can dominate the harmonic analysis  on spherical manifolds. 
The reduction $O(n,R)>S(n+1),\; n>2$ for simplicial manifolds may require 
generalized Casimir operators as exemplified in \cite{KM66} for $n=3$. 
Selection rules for $S(n+1)>C_{n+1}$ eliminate complete representations $D^f$ of the group $S(n+1)$ from the harmonic analysis 
on the sphere $S^{(n-1)}$ when restricted to the simplicial manifold.

To see the  topological variety of the harmonic analysis, compare the tetrahedral 
Platonic 3-manifold ${\cal M}$  analyzed here with the 
dodecahedral Platonic 3-manifold ${\cal M}'$. The homotopy group of Poincare's dodecahedral 3-manifold ${\cal M}'$ is, compare \cite{SE34} pp. 216-8, the binary icosahedral group. It was found in \cite{KR05} that  
the isomorphic group $deck({\cal M}')$ acts exclusively as a subgroup of $SU(2,C)^r$ from the right on the 
sphere $S^3$ in the coordinates  eq. \ref{t33}, with the consequence of a degeneracy 
of the dodecahedral eigenmodes. The multiplicity in the reduction from 
$O(4,R)$ to the subset of  eigenmodes for the dodecahedral 3-manifold is completely resolved by a generalized
Casimir operator. 
The multiplicity analysis in \cite{KR05}, \cite{KR06} shows that 
the lowest  dodecahedral eigenmodes are of degree $(2j)=12$.

Comparison with  the present harmonic analysis for the simplicial 
3-manifold demonstrates a genuine dependence of  the selection rules and the spectrum of eigenmodes  on the 
topology and on the topologically invariant subgroups involved. Corresponding implications 
can be drawn for the use of harmonic analysis in the cosmic topology of 3-space.

\section*{Appendix: Some properties of the group $SU(2,C)$ and its representations.}
The group elements $u \in SU(2,C)$ have the unitary, unimodular and complex conjuation properties
\begin{eqnarray}
 \label{a1}
&&u^{\dagger}=u^{-1},\; {\rm det}(u)=1,
\\ \nonumber
&&\overline{u}=q^{-1}uq,\;
q^{-1}=-q=
\left[
\begin{array}{ll}
 0&1\\
-1&0
\end{array}
\right], \: q^T=q^{-1}
\end{eqnarray}
The irreducible representations of $SU(2,C)$ are the Wigner $D^j$ functions
\begin{equation}
\label{a2}
D^j_{m_1m_2}(u), j=0,1/2,1, \ldots, -j\leq m_i\leq j,\: i=1,2.
\end{equation}
which are homogeneous polynomials with real coefficients of degree $2j$ in the 
complex matrix elements eq. \ref{t33} of $u$. They are explicitly given, \cite{WI59} 
pp. 163-6,
\cite{ED57} pp. 56-67, by
\begin{eqnarray}
 \label{a2a}
&&D^j_{m_1m_2}(z_1,z_2,\overline{z}_1,\overline{z}_2)= 
\left[
\frac{(j+m_1)!(j-m_1)!}{(j+m_2)!(j-m_2)!}\right]^{1/2}
\\ \nonumber
&&\times \sum_{\sigma}\frac{(j+m_2)!(j-m_2)!(-1)^{m_2-m_1+\sigma}} 
{(j+m_1-\sigma)!(m_2-m_1+\sigma)!\sigma!(j-m_2-\sigma)!}
\\ \nonumber
&&\times z_1^{j+m_1-\sigma}\overline{z}_2^{m_2-m_1+\sigma} 
z_2^{\sigma}\overline{z}_1^{j-m_2-\sigma},
\\ \nonumber 
&&j=0,1/2,1,3/2,\ldots
\end{eqnarray}

Unitarity and reality relations of these polynomial representations imply
\begin{eqnarray}
 \label{a3}
&&D^j_{m_1m_2}(u^{-1})= \overline{D^j_{m_2m_1}(u)},\:
\\ \nonumber 
&&\overline{D^j_{m_1m_2}(u)}= D^j_{m_1m_2}(\overline{u}).
\end{eqnarray}
From these equations follows the transposition property 
\begin{equation}
\label{a4}
D^j_{m_1m_2}(u^T)=D^j_{m_2m_1}(u).
\end{equation}
The polynomial homogeneity from eq. \ref{a2a} reads
\begin{equation}
 \label{a5}
D^j_{m_1m_2}(\lambda u)= \lambda^{(2j)}D^j_{m_1m_2}(u).
\end{equation}
To pass from an element $g \in SU(2,C)$ to its character $\chi^j(g)$ 
in the representation $D^j(g)$ we 
first obtain the angle $\phi$ from its trace, 
\begin{equation}
 \label{t54}
\frac{1}{2}\chi^{1/2}(g)= \frac{1}{2} {\rm Trace}(g):= \cos(\phi/2),
\end{equation}
and then for general $j$  use 
\begin{equation}
 \label{t55}
\chi^j(g) =\sum_{m=-j}^{j} \exp(im\phi)
=\frac{\sin(\frac{\phi(2j+1)}{2})}{\sin(\frac{\phi}{2})}.
\end{equation}

\section{Acknowledgement.}
My  thanks are due to  Marcos Moshinsky, UNAM Mexico, who over many years allowed me to  share  his insight into
groups and their representation theory.  Moreover during a visit to Mexico in 2008 he participated in 
the  set-up of the present research project. Many thanks are due to R. Kellerhals, U. Fribourg, 
Switzerland, for 
pointing out  reference \cite{EV04}, and to T. Kramer, U. Regensburg, Germany, 
for doing  most of the algebraic computations for the Tables.

\end{document}